\documentclass[12pt]{article} % use larger type; default would be 10pt

\usepackage[utf8]{inputenc} % set input encoding (not needed with XeLaTeX)

\usepackage{rotating}
\usepackage{geometry} % to change the page dimensions
\usepackage{graphicx} % support the \includegraphics command and options

\usepackage{booktabs} % for much better looking tables
\usepackage{array} % for better arrays (eg matrices) in maths
\usepackage{paralist} % very flexible & customisable lists (eg. enumerate/itemize, etc.)
\usepackage{verbatim} % adds environment for commenting out blocks of text & for better verbatim
\usepackage{subfig} % make it possible to include more than one captioned figure/table in a single float
\usepackage{ mathrsfs, yfonts}
\usepackage{amsfonts}
\usepackage{amssymb}
\usepackage[all,cmtip]{xy}
\usepackage{amsmath, amsthm, amssymb}
\usepackage{amsthm}

\usepackage{fancyhdr} % This should be set AFTER setting up the page geometry
\usepackage{sectsty}
\allsectionsfont{\sffamily\mdseries\upshape} % (See the fntguide.pdf for font help)
\usepackage[nottoc,notlof,notlot]{tocbibind} % Put the bibliography in the ToC
\usepackage[titles,subfigure]{tocloft} % Alter the style of the Table of Contents

\title{Comparison of sheaf cohomology and singular cohomology}
\author{Yehonatan Sella}

\date{}

\begin{document}
\maketitle
\noindent
A fundamental question about sheaf cohomology is how it compares with other cohomology theories. In this paper we focus on the comparison between sheaf cohomology and singular cohomology.

Before recalling the known results, we review some definitions. Let $X$ be a topological space.
\\
\\Definitions: $X$ is \textit{locally contractible} if any open subset $U\subset X$ may be covered by contractible open sets.

The above definition has the following slightly weaker variation: $X$ is \textit{semi-locally contractible} if any open subset $U\subset X$ has an open cover $\{W_i\}$ by open subsets $W_i\subset U$ which are contractible in $U$. That is, the inclusion $W_i\hookrightarrow U$ is homotopic to a constant map. Spanier uses the term ``locally contractible" to denote what we call semi-locally contractible.

The space $X$ is \textit{paracompact} if it is Hausdorff and any cover $\mathcal{U}$ of $X$ has a locally finite refinement, that is, a refinement $\mathcal{U}'$ such that any point $x\in X$ belongs to only finitely many members of $\mathcal{U}'$. For a slightly stronger condition, the space $X$ is \textit{hereditarily paracompact} if every open subset $U\subset X$ is paracompact.
\\
\\We may compare sheaf cohomology and singular cohomology by following a chain of comparisons. Spanier proves that, if $X$ is semi-locally contractible and paracompact and $A$ is any abelian group, then the singular cohomology $H^*(X;A)$ is isomorphic to Alexander cohomology $\bar{H}^*(X;A)$ [1, Corollary 6.9.5]. On the other hand, he proves that if $X$ is paracompact, then Alexander cohomology $\bar{H}^*(X;A)$ is isomorphic to Cech cohomology $\check{H}^*(X;\underline{A})$ [1, Corollary 6.8.8], where $\underline{A}$ denotes the constant sheaf on $X$ with value $A$. Finally, Godement proves that if $X$ is paracompact, then Cech cohomology $\check{H}^*(X;\underline{A})$ is isomorphic to sheaf cohomology $H^*(X;\underline{A})$ [2, Theorem 5.10.1]. Combining these comparisons, we see that if $X$ is semi-locally contractible and paracompact, then singular cohomology $H^*(X;A)$ is isomorphic to sheaf cohomology $H^*(X;\underline{A})$. 
\\
\\The question arises whether all the above conditions on $X$ are necessary. Some kind of condition is certainly necessary, as can be seen by comparing singular cohomology $H^0(X;\mathbb{Z})$, whose rank is the number of path-connected components of $X$, with sheaf cohomology $H^0(X;\underline{Z})$, whose rank is the number of connected components of $X$.

We prove in this paper that we can do away with the paracompactness assumption.
\\
\\Theorem: Let $X$ be a semi-locally contractible topological space. Then for any abelian group $A$ and $n\ge 0$, we have a natural isomorphism between singular cohomology $H^n(X; A)$ and sheaf cohomology $H^n(X; \underline{A})$.
\\
\\This result may not seem new, as it has in fact been stated by some sources, such as Ramanan [3, Theorem 4.14], with the only difference that Ramanan assumed $X$ is locally contractible as opposed to semi-locallly contractibile. However, Ramanan's proof makes implicit use of the assumption that $X$ is hereditarily paracompact.

We first explain Ramanan's method of proof, and the way in which it relies on the assumption that $X$ is hereditarily paracompact. We then modify Ramanan's approach by considering a different flasque resolution of $\underline{A}$ which allows us to remove any paracompactness assumption.

I thank Burt Totaro for pointing out to me the flaw in the original proof and for helpful comments and support. I additionally thank Bhargav Bhatt and Renee Bell for their suggestions.
\section*{Ramanan's method and the reliance on paracompactness}
Let $X$ be a topological space. To fix notation, let $\Delta^n$ denote the standard $n$-simplex, the subspace of $\mathbb{R}^{n+1}$ given by $t_0+\cdots + t_n=1$, $t_i\in \mathbb{R}$. A singular simplex in $X$ (which we may simply refer to as a simplex) is a continuous map $\sigma: \Delta^n\to X$. A face of a singular simplex $\sigma$ is the restriction of $\sigma$ to an $i$-face of $\Delta^n$, for $0\le i\le n$. Let $C_n(X)$ denote the free abelian group on simplices $\sigma:\Delta^n\to X$ and let $C^n(X;A)$ denote the dual $\text{Hom}(C_n(X),A)$.
\\
\\The general method of comparing singular and sheaf cohomology, which encompasses both Ramanan's proof and the approach of this paper, is given by the following lemma.
\\
\\Lemma 0.1: Let $X$ be a semi-locally contractible space. Let $\mathcal{F}^\bullet$ be a complex of sheaves on $X$ with the property that for every open subset $U\subset X$, the complex $\mathcal{F}^\bullet(U)$ is a direct limit of surjections $\mathcal{F}^\bullet_\alpha(U)\to \mathcal{F}^\bullet_\beta(U)$ indexed by $\alpha\le \beta$ in some directed set $D$, equipped with surjective quasi-isomorphisms $\pi_{U,\alpha}:C^\bullet(U;A)\to \mathcal{F}^\bullet_\alpha(U)$ for every $\alpha\in D$ that commute with the direct system. Then $\mathcal{F}^\bullet$ forms a flasque resolution of $\underline{A}$ and the limit map $C^\bullet(X;A)\to \mathcal{F}^\bullet(X)$ is a quasi-isomorphism, inducing the desired isomorphism $H^*(X;A)\cong H^*(X;\underline{A})$.
\\
\\Proof: We first show that the sheaves $\mathcal{F}^n$ are flasque.  Let $U\subset X$ be an open subset. By taking the limit of the surjections $\pi_{U,\alpha}:C^n(U;A)\to \mathcal{F}^n_\alpha(U)$, we obtain a surjection $\pi_U:C^n(U;A)\to \mathcal{F}^n(U)$. We have a commutative diagram 
\[
\xymatrix{
C^n(X;A) \ar[r] \ar[d] & C^n(U;A) \ar[d]
\\
\mathcal{F}^n(X) \ar[r] & \mathcal{F}^n(U)
} 
\]
The top map and the vertical maps are surjections. We conclude the bottom map is a surjection as well, showing that $\tilde{C^n}(-;A)$ is flasque. 

Now, by our assumption that for each $\alpha$, the map $\pi_{U,\alpha}:C^n(U;A)\to \mathcal{F}^n_\alpha(U)$ is a quasi-isomorphism, it follows that the map $\pi_U:C^n(U;A)\to \mathcal{F}^n(U)$ is a quasi-isomorphism as well, since cohomology commutes with direct limits. 

Note that $\underline{A}(U)$ is isomorphic to the kernel of $d:C^0(U;A)\to C^1(U;A)$, hence by the above quasi-isomorphism, $\underline{A}(U)$ is also isomorphic to the kernel of $d:\mathcal{F}^0(U)\to \mathcal{F}^1(U)$. As a convenient notation, let $\hat{C^\bullet}(U;A)$ (resp. $\hat{\mathcal{F}^\bullet}(U)$) denote the result of augmenting $C^\bullet(U;A)$ (resp. $\mathcal{F}^\bullet(A)$) by adding $\underline{A}(U)$ in degree $-1$, with the differential $d_{-1}$ being the natural isomorphism onto the kernel of $d_{0}$. The quasi-isomorphism $\pi_U:C^\bullet(U;A)\to \mathcal{F}^\bullet(U)$ induces a quasi-isomorphism of augmented complexes $\hat{C^\bullet}(U;A)\to \hat{\mathcal{F}^\bullet}(U)$.

Now we show that $\hat{\mathcal{F}^\bullet}$ is acyclic.

Let $U\subset X$ be an open subset. By semi-local contractibility, let $\{U_i\}_{i\in I}$ be a covering of $U$ by open subsets which are contractible in $U$ in the sense that the inclusion $U_i\hookrightarrow U$ is null-homotopic. We then get the following commutative square, for any $i$ and $n$:
\[
\xymatrix{
H^n(\hat{C^\bullet}(U;A)) \ar[r]^{0} \ar[d] & H^n(\hat{C^\bullet}(U_i;A)) \ar[d]
\\
H^n(\hat{\mathcal{F}^\bullet}(U)) \ar[r] & H^n(\hat{\mathcal{F}^\bullet}(U_i))
} 
\]
Since we've shown the vertical maps are isomorphisms, it follows that the bottom map is $0$ as well. That is, any cocycle in $\hat{\mathcal{F}^\bullet}(U)$ restricts to a coboundary in each $\hat{\mathcal{F}^\bullet}(U_i)$, hence being a section of the image sheaf $\text{im}(d_{n-1})(U)$, for $d_{n-1}:\hat{\mathcal{F}^{n-1}}\to \hat{\mathcal{F}^n}$. It follows that the complex of sheaves $\hat{\mathcal{F}^\bullet}$ is acyclic, as desired.
\\
\\Hence, the complex $\mathcal{F}^\bullet$ indeed forms a flasque resolution of $\mathcal{A}$, so $H^n(\mathcal{F}^\bullet(X))$ computes the sheaf cohomology $H^n(X;\underline{A})$. On the other hand, since $\pi_X$ is a quasi-isomorphism, $H^n(C^\bullet(X;A))\cong H^n(\mathcal{F}^\bullet(X))$, which concludes the proof. \qed
\\
\\The natural choice of sheaf $\mathcal{F}^n$ on $X$ admitting a map $C^n(-;A)\to \mathcal{F}^n$ is the sheafification $\tilde{C^n}(-;A)$ of the presheaf $C^n(-;A)$. Ramanan's proof amounts to this proof of lemma 0.1, applied to $\mathcal{F}^n=\tilde{C^n}(-;A)$. Given an open subset $U\subset X$, in order to express $\tilde{C^n}(U;A)$ as a limit of complexes as in lemma 0.1, Ramanan relies on the assumption that the natural map $\pi:C^n(U;A)\to \tilde{C^n}(U;A)$ is surjective. 

Assume for the moment that $\pi$ is surjective. To see how this assumption helps, note that the kernel of $\pi$, which we'll denote $C_0^n(U;A)\subset C^n(U;A)$, consists exactly of cochains which restrict to $0$ over some open cover of $U$. Since $\pi$ is surjective, we identify $\tilde{C^n}(U;A)=C^n(U;A)/C_0^n(U;A)$.

We can reinterpret this as a direct limit as follows. For every open cover $\mathcal{U}$ of $U$, let $C_n^\mathcal{U}(U)$ denote the free abelian group on simplices $\sigma:\Delta^n\to U$ that land in one of the members of $\mathcal{U}$. Let $C^n_\mathcal{U}(U)$ denote the dual group $\text{Hom}(C_n^\mathcal{U}(U),A)$. Then $C^n(U;A)/C_0^n(U;A)$ is precisely the direct limit of the groups $C^n_\mathcal{U}(U)$ over all open covers, ordered by refinement. An argument by repeated barycentric subdivision shows that the inclusion $C_\bullet^\mathcal{U}(U)\to C_\bullet(U)$ is a chain homotopy equivalence [4, Proposition 2.21] for any topological space $U$ and open covering $\mathcal{U}$ of $U$. Thus the dual restriction map $C^\bullet(U;A)\to C^\bullet_\mathcal{U}(U)$ is a chain homotopy equivalence as well, in particular a quasi-isomorphism. So the conditions of lemma 0.1 are met.
\\
\\It remains to justify the surjectivity of the maps $C^n(U;A)\to \tilde{C^n}(U;A)$. If we assume $X$ is hereditarily paracompact, then we can use the following proposition:
\\
\\Proposition 0.2 [3, Proposition 4.12]: Let $X$ be hereditarily paracompact and let $\mathcal{F}$ be a presheaf on $X$. Let $\tilde{\mathcal{F}}$ denote its sheafification. Suppose that the presheaf $\mathcal{F}$ already satisfies gluability in the sense that, for any covering $\mathcal{U}=\{U_i\}$ of an open subset $U\subset X$ and choice of sections $\alpha_i\in \mathcal{F}(U_i)$ which agree on pairwise intersections, there is a section $\alpha\in \mathcal{F}(U)$ which restricts to $\alpha_i$ over $U_i$. Then for any open set $U\subset X$, the natural map $\mathcal{F}(U)\to \tilde{\mathcal{F}}(U)$ is a surjection.
\\
\\In particular, proposition 0.2 may be applied to the presheaf $C^n(-;A)$, which indeed satisfies gluability in the above sense. So if $X$ is hereditarily paracompact, the proof of lemma 0.1 goes through.
\\
\\Ramanan makes use of proposition 0.2 above, which he explicitly states only under the hereditary paracompactness assumption. However, he neglects to carry over that assumption to the statement of the theorem comparing sheaf and singular cohomology.

The question arises: does lemma 0.1 apply to the complex of sheaves $\tilde{C^n}(-;A)$ without the assumption of paracompactness? In particular, if $X$ is locally contractible, are the natural maps $C^n(X;A)\to \tilde{C^n}(X;A)$ guaranteed to be surjections? The answer is no, as the following example illustrates. 
\\
\\Example 0.3: we give an example of a locally contractible topological space $X$ such that the natural map $C^n(X;A)\to \tilde{C^n}(X;A)$ fails to be surjective.
\\
\\Consider $X=\{1,2,3,4,5\}$, with topology given by the basis of open sets $U_1:=\{1,2,3,4\},U_2:=\{2,3,4,5\},\{2,3\},\{3,4\},\{3\}$. This space $X$ is not Hausdorff, but it is locally contractible. Indeed, note that in general, if $Y$ is a topological space and $y\in Y$ is contained in every open subset of $Y$, then $Y$ is contractible to $y$. Since every open subset of $X$ contains $3$, every open subset of $X$ is contractible. 

We now illustrate the failure of surjectivity of the map $C^1(X)\to \tilde{C^1}(X)$. Define $f_1\in C^1(U_1)$ by the following rule on $1$-simplices $\sigma$: if the image of $\sigma$ is contained in $\{2,3\}$ or in $\{3,4\}$, then $f_1(\sigma)=1$. Else $f_1(\sigma)=0$. Similarly, define $f_2 \in C^1(U_2)$ by the following rule on $1$-simplices $\sigma$: if the image of $\sigma$ is contained in $\{2,3\}$ or in $\{3,4\}$, then $f_2(\sigma)=1$. Else $f_2(\sigma)=2$.

Then $f_1$ and $f_2$ have the same germs at all points of the intersection $U_1\cap U_2$, because they agree on the open cover of $U_1\cap U_2$, $\{\{2,3\},\{3,4\}\}$. So the data of $f_1,f_2$ gives a section of $\tilde{C^1}(X)$. However, it doesn't glue to a section $f\in C^1(X)$. If it did, then the requirement that $f$ and $f_2$ have the same germ at $5$ would force $f$ to restrict to $f_2$ on all of $U_2$, since every neighborhood of $5$ contains $U_2$.  Similarly, the requirement that $f$ and $f_1$ have the same germ at $1$ would force $f$ to restrict to $f_1$ on all of $U_1$. But $f_1$ and $f_2$ do not agree on the intersection $U_1\cap U_2$, since there exists a surjective one-simplex $\sigma:I\to U_1\cap U_2$, for which $f_1(\sigma)=0, f_2(\sigma)=1$; we obtain a contradiction. \qed
\\
\\Since surjectivity of the maps $C^n(U;A)\to \tilde{C^n}(U;A)$, a key ingredient of Ramanan's proof, fails, we must find a different argument if we wish to remove the assumption of hereditary paracompactness. This paper modifies Ramanan's method to work under only the assumption that $X$ is semi-locally contractible. In order to do this, we apply lemma 0.1 to a complex of sheaves formed by a finer direct limit which lends itself better to the required generality.
\section*{The proof}
\noindent Let $X$ be a topological space. Let $A$ be an abelian group.
\\
\\Definition: Define a \textit{nesting} on $X$ to be an assignment $\eta$ of an open subset of $X$ to each finite (possibly empty) sequence in $X$, such that:

i) $\eta(\emptyset)=X$.

ii) For any $n\ge 1$ and sequence $x_1,\dots, x_n$ such that $x_1\in \eta(x_2,\dots, x_n)$, the open set $\eta(x_1,\dots, x_n)$ is an open neighborhood of $x_1$ in $X$.

iii) For any nondecreasing function $f\colon \{1,\dots, m\}\to \{1,\dots,n\}$, we have the containment $\eta(x_1,x_2,\dots,x_n)\subset \eta(x_{f(1)},x_{f(2)},\dots,x_{f(m)})$. 
\\
\\Note: a special case of condition iii) is that $\eta(x_1,\dots, x_n)\subset \eta(x_2,\dots, x_n)$ and similarly $\eta(x_1,\dots,x_n)\subset \eta(x_1,\dots, x_{n-1})$ for any sequence $x_1,\dots,x_n$ with $n\ge 1$.
\\
\\Given a simplex $\sigma\colon \Delta^n\to X$, let $b(\sigma)$ denote the image of the barycenter of $\Delta^n$ in $X$. 

Given $n\ge 0$ and a nesting $\eta$ on $X$, let $C_n^{\eta}(X)$ denote the free abelian group on simplices $\sigma\colon \Delta^n\to X$ satisfying the property that for any chain of face inclusions $\sigma_1\subset \cdots \subset \sigma_m\subset \sigma$, the face $\sigma_1$ is contained in $\eta(b(\sigma_1),\dots,b(\sigma_m))$. Define $C^n_{\eta}(X; A)$ to be the dual $\text{Hom}(C_n^{\eta}(X),A)$. Let $\mathcal{C}^n(X;A)$ denote the direct limit of $C^n_{\eta}(X;A)$ over all nestings $\eta$, ordered by refinement.

Note: due to condition iii) of a nesting, in order to show a simplex $\sigma$ lies in $C_n^\eta(X)$, it suffices to check the containment $\sigma_1\subset \eta(b(\sigma_1),\dots,b(\sigma_m))$ only for chains of proper face inclusions $\sigma_1\subsetneq \cdots \subsetneq \sigma_m=\sigma$.
\\
\\Given a map $f\colon X\to Y$ and a nesting $\eta$ on $Y$, define $f^*\eta$ to be the nesting on $X$ given by $f^*\eta(x_1,\dots,x_n)=f^{-1}(\eta(f(x_1),\dots,f(x_n)))$. If $i\colon U\to X$ is the inclusion of an open subset and $\eta$ is a nesting on $X$, define $\eta|_U=i^*\eta$. Note that $C_n^{\eta|_U}(U)\subset C_n^{\eta}(X)$ so we have a natural restriction $C^n_{\eta}(X;A)\to C^n_{\eta|_U}(U;A)$ which induces a map on the direct limits $\mathcal{C}^n(X;A)\to \mathcal{C}^n(U;A)$.

Similarly, if $\eta$ is a nesting on $X$, then given a simplex $\sigma\colon \Delta^n\to X$ that lies in $C_n^\eta(X)$, we have $\partial \sigma\in C_{n-1}^{\eta}(X)$, so that $C_\bullet^\eta(X)$ forms a complex. Dualizing, $C^\bullet_\eta(X)$ forms a complex as well, so by taking direct limits, $\mathcal{C}^\bullet(X)$ forms a complex.
\\
\\
Let $X$ be a semi-locally contractible space. We will apply lemma 0.1 to the complex of presheaves $\mathcal{C}^\bullet(-;A)$ on $X$. For every open subset $U\subset X$, the complex $\mathcal{C}^\bullet(U;A)$ is by definition a direct limit of the complexes $C^\bullet_{\eta}(U;A)$, which come equipped with natural surjections $\pi_U:C^\bullet(U;A)\to C^\bullet_{\eta}(U;A)$.

In step 1, we prove that the presheaf $\mathcal{C}^\bullet(-;A)$ in fact forms a sheaf. Then we need a small-chains argument showing that, for any topological space $U$ and nesting $\eta$ on $U$, the inclusion $C_\bullet^\eta(U)\subset C_\bullet(U)$ is a chain homotopy equivalence, so that the dual restriction $\pi_U: C^\bullet(U;A)\to C^\bullet_\eta(U;A)$ is a chain homotopy equivalence as well, and the conditions of lemma 0.1 are met. This turns out to be significantly more involved than the corresponding argument for $C^{\mathcal{U}}_\bullet(U)\to C_\bullet(U)$, and occupies the bulk of the proof (steps 2 and 3).
\\
\\\textbf{Step 1}: For each topological space $X$, we show that the presheaf $\mathcal{C}^n(-;A)$ on $X$ is in fact a sheaf.

For this it helps to give an equivalent definition of $\mathcal{C}^n(X;A)$. Let $n\ge 0$, $\eta$ a nesting, $x\in X$. Then define $C_{n,x}(X)$ to be the subgroup of $C_n(X)$ consisting of chains all of whose simplices have barycenter $x$, and let $C^{n,x}(X;A)$ be the dual group $\text{Hom}(C_{n,x}(X),A)$. Define $C_{n,x}^{\eta}(X)=C_{n,x}(X)\cap C_n^\eta(X)$ and let $C^{n,x}_{\eta}(X;A)$ be the dual group $\text{Hom}(C_{n,x}^{\eta}(X),A)$. Define $\mathcal{C}^{n,x}(X;A)$ to be the direct of the groups $C^{n,x}_\eta(X;A)$ over all nestings $\eta$ on $X$. In order to think of $\mathcal{C}^{n,x}(-;A)$ as a presheaf on $X$, we define $\mathcal{C}^{n,x}(U;A)=0$ if $x\notin U$.
\\
\\Lemma 1.1: For any open subset $W\subset X$ and $x\in W$, the natural restriction map $\pi\colon \mathcal{C}^{n,x}(X;A)\to \mathcal{C}^{n,x}(W;A)$ is a bijection.
\\
\\Proof: Let $W\subset X$ be an open subset and let $x\in W$. 

The proof of flasqueness in lemma 0.1 can be modified to prove that $\mathcal{C}^{n,x}(-;A)$ is flasque. Thus it remains to show that the restriction $\pi\colon \mathcal{C}^{n,x}(X;A)\to \mathcal{C}^{n,x}(W;A)$ is injective. For this, it suffices to show that, for each nesting $\eta$ on $W$, there is a nesting $\tilde{\eta}$ on $X$, such that $C_{n,x}^{\tilde{\eta}}(X)\subset C_{n,x}^{\eta}(W)$.

We construct the desired nesting $\tilde{\eta}$ on $X$ as follows. Given a sequence $x_1,\dots, x_m$ in $X$, we have two cases. If for some $0\le k\le m$ we have $x_i\in W$ for $i\le k$ and $x_i\notin W$ for $k<i\le m$, then we define $\tilde{\eta}(x_1,\dots,x_m)=\eta(x_1,\dots,x_k)$. Otherwise, define $\tilde{\eta}(x_1,\dots, x_m)=\emptyset$. It is routine to check that $\tilde{\eta}$ is indeed a nesting. To show the containment $C_{n,x}^{\tilde{\eta}}(X)\subset C_{n,x}^{\eta}(W)$, note that if $\sigma$ is a simplex in $C_{n,x}^{\tilde{\eta}}(X)$, then by definition, it follows in particular that $\sigma\subset \tilde{\eta}(b(\sigma))=\tilde{\eta}(x)=\eta(x)$, where the last equality uses the assumption that $x\in W$. Thus $\sigma$ lands in $W$ so for any chain of faces $\sigma_1\subset \cdots \subset \sigma_m=\sigma$, we have $b(\sigma_i)\in W$ for all $i$, so $\sigma_1\subset \tilde{\eta}(b(\sigma_1),\dots, b(\sigma_m))=\eta(b(\sigma_1),\dots, b(\sigma_m))$, showing $\sigma\in C_{n,x}^{\eta}(W)$.
\qed
\\
\\By lemma 1.1, we see that the presheaf $\mathcal{C}^{n,x}(-;A)$ on $X$ is in fact a skyscraper sheaf over $x$, so in particular $\mathcal{C}^{n,x}(X;A)=\mathcal{C}^{n,x}(-;A)_x$, the stalk of $\mathcal{C}^{n,x}(-;A)$ at $x$.
\\
\\Lemma 1.2: The natural product of restriction maps $p\colon \mathcal{C}^n(X;A)\to \Pi_{x\in X}\mathcal{C}^{n,x}(X;A)$ is a bijection.
\\
\\Proof: Again it is easy to see that $p$ is surjective, by considering the square

\[
\xymatrix{
C^n(X;A) \ar[r] \ar[d] & \Pi_{x\in X}C^{n,x}(X;A) \ar[d]
\\
\mathcal{C}^n(X;A) \ar[r] & \Pi_{x\in X}\mathcal{C}^{n,x}(X;A)
} 
\]

It remains to show that $p$ is injective. For this it suffices to show that for any choice of nestings $(\eta^x)_{x\in X}$ on $X$, there exists a single nesting $\eta$ on $X$ such that $C_n^{\eta}(X)\subset \bigoplus_{x\in X} C_{n,x}^{\eta^x}(X)$.

Define $\eta(x_1,\dots,x_m)=\bigcap_{i=1}^m \eta^{x_i}(x_1,\dots, x_i)$. It is routine to verify that $\eta$ is indeed a nesting on $X$. Now suppose $\sigma$ is a simplex in $C_n^\eta(X)$. By definition, for any chain of simplices $\sigma_1\subset \cdots \subset \sigma_m=\sigma$, we have $\sigma_1\subset \eta(b(\sigma_1),\dots,b(\sigma_m))$. By definition of $\eta$, this is contained in $\eta^{b(\sigma_m)}(b(\sigma_1),\dots, b(\sigma_m))$, showing that indeed, $\sigma\in  C_{n,b(\sigma)}^{\eta^{b(\sigma)}}(X)$, as desired. \qed
\\
\\Combining lemmas 1.1 and 1.2, we see that, for any open subset $U\subset X$, we have $\mathcal{C}^{n}(U;A)=\Pi_{x\in U}\mathcal{C}^{n,x}(-;A)_x$, from which we immediately conclude that $\mathcal{C}^n(-;A)$ is a sheaf.
\\
\\\textbf{Step 2}: Let $X$ be a topological space and $\eta$ be a nesting on $X$. Toward our goal of showing that the inclusion map $i\colon C_\bullet^{\eta}(X)\to C_\bullet(X)$ is a chain homotopy equivalence, in this step we prove the preliminary result that, given any chain $\sigma$ in $C_\bullet^{\eta}(X)$ which is a boundary in $C_\bullet(X)$, it is also a boundary in $C_\bullet^{\eta}(X)$.
\\
\\In this step and the next, we will use a notion of simplicial complex that has slightly more structure than the usual one: for us, each face of a simplicial complex will be equipped with an ordering of its vertices, such that the orderings are compatible with restriction to subfaces. This allows us to understand each $k$-face $\sigma$ of $\mathcal{K}$ as a singular simplex $\sigma:\Delta^k\to |\mathcal{K}|$, where $|\mathcal{K}|$ denotes the underlying topological space of $\mathcal{K}$. 

We let $A_n(\mathcal{K})$ denote the free abelian group on the set of $n$-faces of $\mathcal{K}$. By the compatibility of the orderings of the faces, $A_\bullet(\mathcal{K})$ forms a complex with the natural boundary map, and we have an embedding $A_\bullet(\mathcal{K})\to C_\bullet(\mathcal{K})$. We will often implicitly identify $A_\bullet(\mathcal{K})$ as a subcomplex of $C_\bullet(\mathcal{K})$ via this embedding. Define a \textit{facet} of a simplicial complex $\mathcal{K}$ to be a maximal face.
\\
\\The main idea of this step is to prove the following proposition.
\\
\\Proposition 2.1: Let $k\ge 0$ and let $\eta$ be a nesting on $\Delta^k$. Let $i:C_\bullet^\eta(\Delta^k)\to C_\bullet(\Delta^k)$ denote the inclusion. Then there exists a map of complexes $\pi:A_\bullet(\Delta^k)\to C_\bullet^\eta(\Delta^k)$ and a chain homotopy $h:A_n(\Delta^k)\to C_{n+1}(\Delta^k)$ between $i \circ \pi$ and the natural embedding $A_\bullet(\Delta^k)\to C_\bullet(\Delta^k)$, with the property that $h$ maps $A_\bullet(\Delta^k)\cap C^\eta_\bullet(\Delta^k)$ to $C^\eta_\bullet(\Delta^k)$.
\\
\\Once we prove proposition 2.1, our desired result for this step will follow. Indeed, suppose $X$ is a topological space equipped with a compatible nesting $\eta$. Suppose $\sigma\colon \Delta^k\to X$ is a simplex such that $\partial \sigma \in C_\bullet^\eta(X)$. In other words, we have $\partial \Delta^k\in C_\bullet^{\sigma^*\eta}(\Delta^k)$. Let $\pi\colon A_\bullet(\Delta^k)\to C_\bullet(\Delta^k)$ and $h\colon A_n(\Delta^k)\to C_{n+1}(\Delta^k)$ be as in Proposition 2.1, with respect to the nesting $\sigma^*\eta$ on $\Delta^k$. Then $\partial(\pi(\Delta^k)+h(\partial \Delta^k)) = \pi(\partial \Delta^k)+(\partial \Delta^k - \pi(\partial \Delta^k))=\partial \Delta^k$. We know that $\pi(\Delta^k)\in C_\bullet^{\sigma^*\eta}(\Delta^k)$ and also $h(\partial \Delta^k)\in C_\bullet^{\sigma^*\eta}(\Delta^k)$ since $\partial \Delta^k\in C_\bullet^{\sigma^*\eta}(\Delta^k)$. So $\partial \Delta^k$ is indeed a boundary of a chain in $C_\bullet^{\sigma^*\eta}(\Delta^k)$. Pushing forward by $\sigma$, we see that $\partial \sigma$ is the boundary of a chain in $C_\bullet^\eta(X)$.
\\
\\We now work toward the proof of proposition 2.1. Roughly speaking, the map $\pi$ will consist of first performing barycentric subdivision on $\Delta^k$ to break it up into sufficiently small subsimplices, and then deforming these subsimplices so that their barycenters satisfy the needed conditions.
\\
\\We mention a few operations on simplices and simplicial complexes which we will use in this step. The following constructions 2.2-2.4 are taken from Hatcher [4, Section 2.1], with slight modification.
\\
\\Given a convex space $X$, and points $v_0,\dots, v_k\in X$, let $[v_0,\dots,v_k]\subset X$ denote the linear ordered simplex in $X$ with vertices $v_0,\dots, v_k$ in that order. Given a point $x\in X$ and a linear ordered simplex $\sigma=[v_0,\dots,v_k]\subset X$, we let $C_x(\sigma)$ denote the cone on $\sigma$ with vertex $x$, that is, the linear ordered simplex $[x,v_0,\dots, v_k]$.
\\
\\Let $\mathcal{K}$ be denote a simplicial complex.
\\
\\Construction 2.2: We construct the barycentric subdivision $S(\mathcal{K})$, a simplicial complex structure on $|\mathcal{K}|$. We also obtain a corresponding map of complexes $S: A_\bullet(\mathcal{K})\to A_\bullet(S(\mathcal{K}))$.
\\
\\Proof: It suffices to define the simplicial complexes $S(\Delta^k)$ for $k\ge 0$ and to show the compatibility condition that $S(\Delta^k)$ restricts to $S(\tau)$ over a face $\tau$ of $\Delta^k$. Once we do this, then we define $S(\mathcal{K})$ in general as the union of $S(\sigma)$ over the faces $\sigma$ of $\mathcal{K}$. Furthermore, given a face $\alpha$ of $\mathcal{K}$, we define the chain $S(\alpha)\in A_\bullet(S(\mathcal{K}))$ as $\sum_{\sigma} \pm 1 \cdot \sigma$, where $\sigma$ ranges over the facets of the simplicial complex $S(\alpha)$ and the $\pm 1$ is the orientation of $\sigma$ relative to $\alpha$. The compatibility condition on faces guarantees that this a map of complexes.

We proceed to construct $S(\Delta^k)$. Define $S(\Delta^0)=\Delta^0$. For $k>0$, define $S(\Delta^k)$ inductively as the union of the cones $C_b(\sigma)$, where $b$ is the barycenter of $\Delta^k$ and $\sigma$ ranges over all facets of $S(\partial \Delta^k)$. By induction we may assume the faces $\tau$ of $S(\partial \Delta^k)$ have compatible orderings. Thus the cones $C_b(\sigma)$ have compatible orderings as well, so $S(\Delta^k)$ forms a simplicial complex. It also follows by induction that $S(\Delta^k)$ restricts to $S(\tau)$ over a face $\tau$, as desired. \qed
\\
\\Let $S^n$ denote $n$-fold barycentric subdivision $S^{\circ n}$.
\\
\\We now define two simplicial complex structures on $|\mathcal{K}|\times I$, accompanied by corresponding chain homotopies. As in construction 2.1, we need only define each construction on simplices and verify the compatibility of the construction with respect to faces.
\\
\\Construction 2.3: We construct a simplicial complex structure $T_n(\mathcal{K})$ on $|\mathcal{K}|\times I$ which restricts to $S^n(\mathcal{K})$ over $|\mathcal{K}|\times \{0\}$ and to $\mathcal{K}$ over $|\mathcal{K}|\times \{1\}$. We obtain corresponding maps on chains $T_n:A_k(\mathcal{K})\to A_{k+1}(T_n(\mathcal{K}))$, which form a chain homotopy between $(i_0)_*\circ S^n$ and $(i_1)_*$, where $i_0$ and $i_1$ denote the inclusions $S^n(\mathcal{K})\to T_n(\mathcal{K})$ at $0$ and $\mathcal{K}\to T_n(\mathcal{K})$ at $1$, respectively.
\\
\\Proof: We first construct $T:=T_1$. Define $T(\Delta^k)$ recursively by letting $T(\Delta^0)= \Delta^1\cong \Delta^0\times I$, and if $k>0$, letting $T(\Delta^k)$ be the union of the cones $C_b(\sigma)$, where $b$ is the barycenter of $\Delta^k\times \{0\}$ and $\sigma$ is a facet of $\Delta^k \times \{1\} \cup T(\partial \Delta^k)$. Again by an inductive argument, the $C_b(\sigma)$ have compatible orderings, and $T(\Delta^k)$ restricts to $T(\tau)$ over $|\tau|\times I$.

We can use $T$ to define the $T_n$ for $n>1$. Indeed, let $T_n(\mathcal{K})$ to be result of gluing the simplicial complexes $T(S^{n-i}(\mathcal{K}))$, ranging over $i=1,\dots, n$, where we rescale the intervals to view $T(S^{n-i}(\mathcal{K}))$ as a simplicial complex structure on $|\mathcal{K}|\times [\frac{i-1}{n},\frac{i}{n}]$. \qed
\\
\\Construction 2.4: We define a simplicial complex structure $P(\mathcal{K})$ on $|\mathcal{K}|\times I$ whose restrictions to $|\mathcal{K}|\times\{0\}$ and $|\mathcal{K}|\times \{1\}$ are both $\mathcal{K}$. The corresponding maps of chains $P: A_k(\mathcal{K})\to C_{k+1}(P(\mathcal{K}))$ form a chain homotopy between $(i_0)_*$ and $(i_1)_*$, where here $i_0$ and $i_1$ represent the inclusions of $\mathcal{K}$ in $P(\mathcal{K})$ at $0$ and $1$, respectively.
\\
\\Proof: Define $P(\Delta^k)$ to be the union of simplices $[v_0,\dots,v_i,w_i,\dots,w_k]$ for $i=0,\dots,k$, where $v_0,\dots,v_k$ are the ordered vertices of $\Delta^k\times \{0\}$ and $w_0,\dots, w_k$ are the ordered vertices of $\Delta^k\times \{1\}$. It follows directly from construction that the facets of $P(\Delta^k)$ have compatible orderings and that $P(\Delta^k)$ restricts to $P(\tau)$ over $|\tau|\times I$ for any face $\tau$ of $\Delta^k$. \qed
\\
\\We now study conditions under which a singular chain in $|\mathcal{K}|$ can be ``deformed" to a chain that lies in $C_\bullet^\eta(|\mathcal{K}|)$.
\\
\\Definition: Let $S$ be a set of faces of $\mathcal{K}$. Then a \textit{compatible }$\mathit{\eta}$\textit{-covering} of $S$ is a choice, for each $i$-face $\alpha\in S$, of a contractible subset $W(\alpha)\subset |\mathcal{K}|$ containing $\alpha$, and an element $t(\alpha) \in W(\alpha)$ such that:

i)If $\alpha\in S$ is a $0$-face, then $t(\alpha)=\alpha$.

ii) Whenever $\beta\subset \alpha$ is a containment of faces in $S$, we have $W(\beta)\subset W(\alpha)$.

iii) For any chain of faces $\alpha_1\subset \cdots \subset \alpha_k$  in $S$, we have $W(\alpha_1)\subset \eta(t(\alpha_1),\dots,t(\alpha_k))$.
\\
\\If $|\mathcal{K}|$ is embedded in an ambient space $X$, we will more generally allow $W(\alpha)\subset X$. Define a compatible $\eta$-covering of $\mathcal{K}$ to be a compatible $\eta$-covering of the set of all faces of $\mathcal{K}$.
\\
\\We remark that if $\mathcal{K}$ is a simplicial complex which is a union of subsimplicial complexes $\mathcal{K}_1$ and $\mathcal{K}_2$, then given compatible $\eta$-coverings $(W_1,t_1)$ of $\mathcal{K}_1$ and $(W_2,t_2)$ of $\mathcal{K}_2$ which agree along the intersection $\mathcal{K}_1\cap \mathcal{K}_2$, they may be glued to a compatible $\eta$-covering of $\mathcal{K}$. Indeed, this uses the fact that any chain of faces in $\mathcal{K}$ is contained in a facet of $\mathcal{K}$, which must either reside completely in $\mathcal{K}_1$ or in $\mathcal{K}_2$.
\\
\\The next result explains our interest in compatible $\eta$-coverings.
\\
\\Construction 2.5: Given a compatible $\eta$-covering $(W,t)$ of $\mathcal{K}$, we construct a map of complexes $\delta\colon A_\bullet(\mathcal{K})\to C_\bullet^{\eta}(|\mathcal{K}|)$ such that for any face $\alpha$ of $\mathcal{K}$, $\delta(\alpha)$ is a simplex contained in $W(\alpha)$ with barycenter $t(\alpha)$. Furthermore, if $\alpha$ is a face of $\mathcal{K}$ with the property that $t(\beta)=b(\beta)$ for every face $\beta$ of $\alpha$, then $\delta(\alpha)=\alpha$.
\\
\\Note: the fact that $\delta$ has image in $C_\bullet^\eta(\mathcal{K})$ is automatic from condition iii) of a compatible $\eta$-covering and the requirement that $\delta(\alpha)$ is contained in $W(\alpha)$ with barycenter $t(\alpha)$.
\\
\\Proof: Slightly rephrasing, we construct a continuous function $f\colon |\mathcal{K}|\to |\mathcal{K}|$ such that for every face $\alpha$ of $\mathcal{K}$, $f(\alpha)\subset W(\alpha)$ and $f(b(\alpha))=t(\alpha)$. Furthermore, if $\alpha$ is an $i$-face of $\mathcal{K}$ with the property that $t(\beta)=b(\beta)$ for every face $\beta$ of $\alpha$, then $f|_{\alpha}=\text{id}$. Once we define such a function $f$, we can recover the map of complexes $\delta$ by $\delta(\alpha)=f \circ \alpha:\Delta^i\to |\mathcal{K}|$.

We define $f$ recursively on the $n$-skeletons of $\mathcal{K}$. Define $f(x)=x$ for $x$ in the $0$-skeleton of $|\mathcal{K}|$. By condition i) of a compatible $\eta$-covering, indeed $f(b(x))=t(x)$. Now let $n>0$ and suppose $f$ is defined on the $(n-1)$-skeleton of $|\mathcal{K}|$. We wish to extend $f$ to the $n$-skeleton of $|\mathcal{K}|$. It suffices to extend $f$ from $\partial \alpha$ to $\alpha$ for any given $n$-face $\alpha$ of $\mathcal{K}$.

If every face $\beta$ of $\alpha$ has barycenter $t(\beta)$, then we can assume by induction that $f$ is the identity on $\partial \alpha$, and define $f$ to be the identity on $\alpha$.

Now assume the simplex $\alpha$ is not of that form. We note by induction and by condition ii) of a compatible $\eta$-covering that $f(\partial \alpha)\subset W(\alpha)$. By contractibility of $W:= W(\alpha)$, fix a homotopy $H\colon W\times I\to W$ between the identity and the constant map $t(\alpha)$. Then we extend $f$ from $\partial \alpha$ to $\alpha$ by mapping the line from $x$ to $b(\alpha)$ to the path $H_{f(x)}\colon I\to W$ from $f(x)$ to $t(\alpha)$, for any $x\in \partial \alpha$. Then $f(\alpha)\subset W(\alpha)$ and $f(b(\alpha))=t(\alpha)$, as desired. \qed
\\
\\Of course not every simplicial complex admits a compatible $\eta$-covering, but the next lemma will show that it will upon sufficiently many applications of barycentric subdivision. This will make use of Lebesgue's number lemma [5, Lemma 27.5], which states that for any compact space $X$ and covering $\mathcal{U}$ of $X$, there exists a number $\delta>0$ such that any subset of $X$ with diameter at most $\delta$ is contained in some member of $\mathcal{U}$. We will apply this with the understanding that the diameters of the facets in the barycentric subdivision $S^n(\mathcal{K})$ uniformly converge to zero.
\\
\\First we introduce some terminology. Say a set $S$ of faces of $\mathcal{K}$ is upwards closed if, for any inclusion of faces $\beta\subset \alpha$ of $\mathcal{K}$ where $\beta\in S$, then $\alpha\in S$.

We also introduce the following notation: given $n\ge 0$ and a face $\alpha$ of $S^n(\mathcal{K})$, let $F^{\mathcal{K}}(\alpha)$ denote the smallest face of $\mathcal{K}$ containing $\alpha$.
\\
\\Lemma 2.6: a) There exists an $n$ such that $S^n(\mathcal{K})$ admits a compatible $\eta$-covering $(W,t)$.

b) Furthermore, if $S$ is an upwards closed set of faces of $\mathcal{K}$ and we are given a compatible $\eta$-covering $(W_0,t_0)$ of $S$, then we can pick $(W,t)$ in a) so that for any face $\alpha$ of $S^n(\mathcal{K})$ for which $F^\mathcal{K}(\alpha)\in S$, we have $W(\alpha)=W_0(F^\mathcal{K}(\alpha))$ and $t(\alpha)=t_0(F^\mathcal{K}(\alpha))$.
\\
\\Proof: We use the following inductive step. Suppose we are given $n\ge 0$ and a compatible $\eta$-covering $(W_0,t_0)$ of an upwards closed subset $S$ of faces of $\mathcal{K}$. Then let $S'$ be the set of faces of $\mathcal{K}$ which are maximal faces out of those not in $S$. Given a number $n\ge 0$, let $\bar{S}$ (suppressing the dependency on $n$) be the set of faces $\alpha$ of $S^{n}(\mathcal{K})$ such that either $\alpha$ is a facet of $S^{n}(\tau)$ for some $\tau\in S'$ or $F^\mathcal{K}(\tau)\in S$. Then for some $n\ge 0$, there is a compatible $\eta$-covering $(W,t)$ of the set $\bar{S}$. Furthermore, if $F^{\mathcal{K}}(\alpha)\in S$, then $W(\alpha)=W_0(F^\mathcal{K}(\alpha)), t(\alpha)=t_0(F^\mathcal{K}(\alpha))$.

Note that the top dimension of a face not in $S$ goes down after applying the inductive step and replacing $\mathcal{K}$ with $S^n(\mathcal{K})$ and $S$ with $\bar{S}$, so this gives rise to a terminating process. Once we prove the inductive step, a) will follow by beginning with $S=\emptyset$, and b) will follow by beginning with the given partial compatible $\eta$-covering $(W_0,t_0)$ on the given set $S$.
\\
\\We now prove the inductive step. Let $\tau\in S'$. By definition of $S'$, note that for any strict chain $\tau = \tau_1\subsetneq \tau_2 \subsetneq \cdots \subsetneq \tau_k$ of faces of $\mathcal{K}$, each $\tau_i$ is in $S$ for $i>1$ and we have $\tau\subset \tau_2\subset W_0(\tau_2)\subset \eta(t_0(\tau_2),\dots, t_0(\tau_k))$, since $(W_0,t_0)$ is a compatible $\eta$-cover of $S$. So for any $t\in \tau$, the set $\eta(t,t_0(\tau_2),\dots,t_0(\tau_k))$ is an open neighborhood of $t$.

Now for each $t\in \tau$, let $U_t$ be the intersection of the face $\tau$ with the open neighborhoods $\eta(t,t_0(\tau_2),\dots, t_0(\tau_k))$, over all chains $\tau = \tau_1\subsetneq \tau_2\subsetneq \cdots \subset \tau_k$. There are finitely many such chains, so $U_t$ is open in $\tau$. Note that by condition iii) of a nesting, $U_t$ is in fact contained in $\eta(t,t_0(\tau_2),\dots,t_0(\tau_k))$ for every non-strict chain $\tau=\tau_1\subset \tau_2 \subset \cdots \subset \tau_k$ as well. Now define $W_t$ to be a contractible neighborhood of $t$ in $\tau$ contained in $U_t$. Then by compactness of $\tau$ and Lebesgue's number lemma, there exists an $n$ such that every facet of $S^{n}(\tau)$ is contained in $W_t$ for some $t\in \tau$. We choose $n$ uniformly so that this works for all $\tau$ in $S'$.

We now define the compatible $\eta$-covering $(W,t)$ on the set $\bar{S}$. Let $\alpha$ be a face in $\bar{S}$. 

Case i): $\alpha$ is a facet of $S^n(\tau)$ for some $\tau\in S'$. Then by choice of $n$, we can fix $t(\alpha)\in \tau$ such that $\alpha$ is contained in $W_{t(\alpha)}$. If $\alpha$ is a $0$-face, we can choose $t(\alpha)=\alpha$. We define $W(\alpha)=W_{t(\alpha)}$.

Case ii): $F^\mathcal{K}(\alpha)\in S$. Then define $W(\alpha)=W_0(F^\mathcal{K}(\alpha)), t(\alpha)=t_0(F^\mathcal{K}(\alpha))$.
\\
\\We prove that $(W,t)$ is a compatible $\eta$-covering of $\bar{S}$. By construction, $W(\alpha)$ contains $\alpha$ and $t(\alpha)$, and if $\alpha$ is a $0$-face, $t(\alpha)=\alpha$. 

To prove condition ii) of a compatible $\eta$-covering, let $\beta\subsetneq \alpha$ be faces in $\bar{S}$. We show $W(\beta)\subset W(\alpha)$.

Case 1: $\beta$ is a facet of $S^n(\tau)$ for some $\tau\in S'$. If $\alpha$ were a facet of $S^n(\sigma)$ for some $\sigma\in S'$, this would imply $\tau\subsetneq \sigma$, contradicting the definition of $S'$. Hence $F^\mathcal{K}(\alpha)\in S$. Then
\begin{equation*}
W(\beta)=W_{t(\beta)}\subset \tau=F^\mathcal{K}(\beta)\subset F^\mathcal{K}(\alpha)\subset W_0(F^\mathcal{K}(\alpha))=W(\alpha)
\end{equation*}
Case 2: $F^\mathcal{K}(\beta), F^\mathcal{K}(\alpha)\in S$. Then
\begin{equation*}
W(\beta)=W_0(F^\mathcal{K}(\beta))\subset W_0(F^\mathcal{K}(\alpha))=W(\alpha)
\end{equation*}
since $W_0$ is a compatible covering of $S$ and $F^\mathcal{K}(\beta)\subset F^\mathcal{K}(\alpha)$.
\\
\\Now we prove condition iii) of a compatible $\eta$-covering.

Let $\alpha_1\subsetneq \alpha_2 \subsetneq \cdots \subsetneq \alpha_k$ be a strict chain of faces in $\bar{S}$. Again we have two cases.

Case 1: $\alpha_1$ is a facet of $S^n(\tau)$ for some face $\tau\in S'$. It follows that $F^\mathcal{K}(\alpha_i)\in S$ for $i>1$. Then 
\begin{equation*}
W(\alpha_1)=W_{t(\alpha_1)}\subset  \eta(t(\alpha_1),t_0(F^\mathcal{K}(\alpha_2)),\dots, t_0(F^\mathcal{K}(\alpha_k))=\eta(t(\alpha_1),t(\alpha_2),\dots,t(\alpha_k))
\end{equation*}
where we use the construction of the sets $W_t$ and the fact that we have a chain $\tau\subset F^\mathcal{K}(\alpha_2)\subset \cdots \subset F^\mathcal{K}(\alpha_k)$.

Case 2: $F^{\mathcal{K}}(\alpha_i)\in S$ for all $i$. Then
\begin{equation*}
W(\alpha_1)=W_0(F^\mathcal{K}(\alpha_1))\subset \eta(t_0(F^\mathcal{K}(\alpha_1)),\dots, t_0(F^\mathcal{K}(\alpha_k))) = \eta(t(\alpha_1),\dots, t(\alpha_k))
\end{equation*}
where we use the fact that $(W_0,t_0)$ is a compatible $\eta$-covering, and the fact that we have a chain $F^\mathcal{K}(\alpha_1)\subset \cdots \subset F^\mathcal{K}(\alpha_k)$.
\qed
\\
\\We have the following construction as a corollary.
\\
\\Corollary 2.7: We construct a map of complexes $\pi\colon A_\bullet(\mathcal{K})\to C_\bullet^{\eta}(|\mathcal{K}|)$ which is the identity on $0$-simplices.
\\
\\Proof: By lemma 2.6a), fix $n$ large enough so that $S^n(\mathcal{K})$ admits a compatible $\eta$-covering $(W,t)$. Let $\delta\colon A_\bullet(S^n(\mathcal{K}))\to C_\bullet^{\eta}(|\mathcal{K}|)$ be the map of complexes given by construction 2.5 with respect to $(W,t)$. Let $\pi=\delta \circ S^n$. By construction, $\delta$ is the identity on $0$-simplices, so $\pi$ is as well. \qed
\\
\\We now concentrate on the case of a simplex $\Delta^k$, equipped with a nesting $\eta$. Corollary 2.7 gives us the map $\pi$ needed for proposition 2.1. In order to obtain the homotopy $h$, we will need to refine our approach by considering the following mapping-cone construction.

Let $\mathcal{K}\subset \Delta^k$ be the sub-simplicial complex consisting of all faces of $\Delta^k$ which are in $C_\bullet^\eta(\Delta^k)$. Recall from constructions 2.3 and 2.4 the simplicial complex structures $P(\mathcal{K})$ and $T_n(\mathcal{K})$ on $|\mathcal{K}|\times I$. For our purposes, we regard $P(\mathcal{K})$ as a simplicial complex structure on $|\mathcal{K}|\times [0,1]$ and $T_n(\mathcal{K})$ as a simplicial complex structure on $|\mathcal{K}|\times [1,2]$. Let $\mathcal{L}$ be the result of gluing $\Delta^k$ to $P(\mathcal{K})$ by identifying $\mathcal{K}\times \{0\}$ with $\mathcal{K}\subset \Delta^k$. Now, given $n\ge 0$, let $\mathcal{L}'_n$ be the result of further gluing $T_n(\mathcal{K})$ to $S^n(\mathcal{L})$ along $|\mathcal{K}|\times \{1\}$, as both simplicial complexes restrict to $S^n(\mathcal{K})$ over $|\mathcal{K}|\times \{1\}$. 

Let $p:|\mathcal{K}|\times [0,2]\to |\mathcal{K}|$ be the projection, and let $q:|\mathcal{L}'_n|\to \Delta^k$ be the map induced by $\text{id}\sqcup p: \Delta^k \sqcup |\mathcal{K}|\times [0,2]\to \Delta^k$. Equip $|\mathcal{L}'_n|$ with the nesting $q^*\eta$. The simplicial complex $\mathcal{L}'_n$ serves as a way of carrying out barycentric subdivision while also keeping track of a homotopy between $S^n(\mathcal{K})$ and $\mathcal{K}$. The prism $S^n(P(\mathcal{K}))$ serves as a kind of buffer between $S^n(\Delta^k)$ and $T_n(\mathcal{K})$, giving us flexibility to deform the former while not affecting the latter. We now emulate lemma 2.6 with respect to $\mathcal{L}'_n$ instead of $S^n(\Delta^k)$.
\\
\\Lemma 2.8: There exists an $n$ such that $\mathcal{L}'_n$ admits a compatible $q^*\eta$-covering $(W,t)$. Furthermore, we can choose $(W,t)$ such that, for any face $\alpha$ of $\mathcal{K}\times \{2\}\subset \mathcal{L}'_n$, we have $t(\alpha)=b(\alpha)$.
\\
\\Proof: We begin by focusing on $\mathcal{L}$. Let $S$ be the set of all faces of $P(\mathcal{K})$ which intersect $\mathcal{K}\times \{1\}$. So $S$ is upwards-closed in $\mathcal{L}$. We let $(W_0,t_0)$ be the covering of $S$ given by $W_0(\alpha)=p_*(\alpha)\times [0,2], t(\alpha)=b(\alpha)$. It is routine to check that $(W_0,t_0)$ is indeed a compatible $q^*\eta$-covering of $S$ using the fact that, by definition, faces of $\mathcal{K}$ are in $C_\bullet^\eta(\Delta^k)$. Thus by lemma 2.6b), there is a number $n\ge 0$ and a compatible $q^*\eta$-covering $(W_1,t_1)$ of  $S^n(\mathcal{L})$ such that for any $\alpha\in S^n(\mathcal{L})$ with $F^{\mathcal{L}}(\alpha)\in S$, we have $W_1(\alpha)=W_0(F^{\mathcal{L}}(\alpha))=F^{\mathcal{K}}(p_*\alpha)\times [0,2]$ and $ t_1(\alpha)=t_0(F^{\mathcal{L}}(\alpha))=b(F^\mathcal{L}(\alpha))$. In particular, this holds whenever $\alpha$ is a face of $S^n(\mathcal{K})\times \{1\}$, in which case we may simplify $t_1(\alpha)=(b(F^{\mathcal{K}}(p_*\alpha),1)\in |\mathcal{K}|\times [0,2]$.

On the other hand, we define a compatible $q^*\eta$-covering $(W_2,t_2)$ of $T_n(\mathcal{K})\subset \mathcal{L}'_n$ by $W_2(\alpha)=F^\mathcal{K}(p_*\alpha)\times [0,2], t_2(\alpha)=(b(F^\mathcal{K}(p_*\alpha)),\pi_2(b(\alpha)))\in |\mathcal{K}|\times [0,2]$, where $\pi_2:\mathcal{K}\times [0,2]\to [0,2]$ denotes the projection. We note that the covering $(W_1,t_1)$ agrees with the covering $(W_2,t_2)$ over $S^n(\mathcal{K})\times {1}$. Thus we may glue them together to obtain a single compatible $q^*$-eta covering $(W,t)$ of $\mathcal{L}'_n$. Furthermore, for any face $\alpha$ of $\mathcal{K}\times\{2\}\subset \mathcal{L}'_n$, we have $t(\alpha)=t_2(\alpha)=(b(F^\mathcal{K}(p_*\alpha)),2)=(b(p_*\alpha),2)=b(\alpha)$, as desired. \qed
\\
\\As a corollary, we may refine corollary 2.7 in a step toward proposition 2.1.
\\
\\Corollary 2.9: There exists a map of complexes $\pi: A_\bullet(\Delta^k)\to C_\bullet^\eta(\Delta^k)$ and a homotopy $h_0: A_n(\mathcal{K})\to C_{n+1}^\eta(\Delta^k)$ between $\pi|_{A_\bullet(\mathcal{K})}$ and the natural embedding $A_\bullet(\mathcal{K})\to C_\bullet^\eta(\Delta^k)$.
\\
\\Proof: Let $\delta: A_\bullet(\mathcal{L}'_n)\to C_\bullet^{q^*\eta}(|\mathcal{L}'_n|)$ be the map given by construction 2.5 with respect to the covering $(W,t)$ of $\mathcal{L}'_n$ guaranteed by lemma 2.8. We note that $\delta$ is the identity restricted to $\mathcal{K}\times \{2\}$ because lemma 2.8 guarantees $t(\alpha)=b(\alpha)$ for all faces $\alpha$ of $\mathcal{K}\times \{2\}$, and construction 2.5 ensures that faces with this property are preserved by $\delta$. 

Let $\delta'=q_*\circ \delta:A_\bullet(\mathcal{L}'_n)\to C_\bullet^\eta(\Delta^k)$. As in corollary 2.7, we obtain $\pi$ by $\pi=\delta'\circ S^n$.

Let $h'= S^n\circ P + T_n\colon A_k(\mathcal{K})\to A_{k+1}(\mathcal{L}'_n)$. Then $h'$ is a chain homotopy between the maps $S^n \circ (i_0)_*\colon A_\bullet(\mathcal{K})\to A_\bullet(\mathcal{L}'_n)$ and $(i_2)_*\colon A_\bullet(\mathcal{K})\to A_\bullet(\mathcal{L}'_n)$. Applying $\delta'$, we obtain a chain homotopy $h_0:=\delta'\circ h'$ between the maps $\delta'\circ S^n \circ (i_0)_*$ and $\delta'\circ (i_2)_*$. But $\delta'\circ S^n \circ (i_0)_*$ is exactly $\pi|_{A_\bullet(\mathcal{K})}$ by definition, and $\delta'\circ (i_2)_*:A_\bullet(\mathcal{K})\to C_\bullet^\eta(\mathcal{K})$ is the natural embedding since $\delta$ is the identity on $\mathcal{K}\times \{2\}$. \qed
\\
\\Finally, we can extend $h_0$ to the desired homotopy $h\colon A_n(\Delta^k)\to C_{n+1}(\Delta^k)$ between $i \circ \pi$ and the natural embedding. We do this as follows. If $\tau$ is a face of $\mathcal{K}$, define $h(\tau)=h_0(\tau)$. Given an $i$-face $\tau$ of $\Delta^k$ which is not a face of $\mathcal{K}$, we can assume by induction that for $\alpha$ a $j$-face with $j<i$, the chain homotopy equation $\partial h(\alpha)=\alpha-\pi(\alpha)-h(\partial \alpha)$ holds. Thus $\tau-\pi(\tau)- h(\partial \tau)$ is a cycle. By contractibility of $\Delta^k$, let $h(\tau)$ be some chain in $\Delta^k$ whose boundary is $\tau-\pi(\tau)- h(\partial \tau)$. 

Since $h$ restricts to $i\circ h_0$ over $A_\bullet(\mathcal{K})$, this proves proposition 2.1.
\\
\\\textbf{Step 3} We now build upon step 2 to prove the stronger result that for any space $X$ equipped with a compatible nesting $\eta$, the inclusion $i:C_\bullet^{\eta}(X)\to C_\bullet(X)$ is a chain homotopy equivalence. 

We first focus on simplices $\Delta^k$. For each pair $(\Delta^k,\eta)$, where $k\in \mathbb{N}$ and $\eta$ is a nesting on $\Delta^k$, we fix a map of complexes $\pi^{(\Delta^k,\eta)}\colon A_\bullet(\Delta^k)\to C_\bullet^\eta(\Delta^k)$ as given by corollary 2.7. If $\eta$ is understood, we will suppress it from the notation and write $\pi^{\Delta^k}$.
\\
\\For flexibility in notation we extend the $\pi$ construction to any simplex, that is, any convex hull $\tau$ of affinely independent points in Euclidean space, equipped with an ordering of the vertices. Let $f:\Delta^i\cong \tau$ be the order-preserving linear isomorphism. Then define $\pi^{(\tau,\eta)}:A_\bullet(\tau)\to C^\eta_\bullet(|\tau|)$ to be the map given by conjugating the map of complexes $\pi^{(\Delta^i,f^*\eta)}: A_\bullet(\Delta^i)\to C_\bullet^{f^*\eta}(\Delta^i)$ by the isomorphism $f$.
\\ 
\\The maps of complexes $\pi^{\tau}$ may not be compatible, in the sense that, if $\tau$ is a face of $\sigma$, $\pi^{\sigma}|_{A_\bullet(\tau)}$ may not equal $\pi^\tau$. It turns out we can get around this issue by a purely formal procedure, given the result of step 2.
\\
\\Construction 3.1: We construct for each proper inclusion of simplices $\tau\subsetneq \sigma$ a chain homotopy $\epsilon^\sigma_\tau:A_n(\tau)\to C_{n+1}^\eta(|\sigma|)$ between the maps $\pi^\tau$ and $\pi^\sigma|_{A_\bullet(\tau)}$, satisfying:

i) (Naturality on isomorphisms) Given a chain of face inclusions $\tau\subsetneq \sigma$ and an order-preserving linear isomorphism $f\colon \Delta^i\cong \sigma$, the homotopy $\epsilon^{\sigma}_{\tau}$ is given by conjugating $\epsilon^{\Delta^i}_{f^{-1}(\tau)}$ by $f$.

ii) (Additivity) For any chain of face inclusions $\alpha\subsetneq \tau\subsetneq \sigma$, we have $\epsilon^\sigma_\alpha=\epsilon^\sigma_\tau|_{A_\bullet(\alpha)}+\epsilon^\tau_\alpha$.
\\
\\Proof: Condition i) serves as a definition in the case that $\sigma$ is not a standard simplex. So we may assume $\sigma=\Delta^k$.

Assume that we have defined $\epsilon^{\sigma'}_{\tau'}$ whenever the combined dimension of $\tau'$ and $\sigma'$ is smaller than that of $\sigma$ and $\tau$. Then we use condition ii) as a definition of $\epsilon^\sigma_\tau$ on proper faces $\alpha$ of $\tau$, by the recursive formula
\begin{equation*}
\epsilon^\sigma_\tau(\alpha)=\epsilon^\sigma_\alpha(\alpha)-\epsilon^\tau_\alpha(\alpha)
\end{equation*}
We have not yet defined $\epsilon^\sigma_\tau$ at $\tau$, but for any proper face $\alpha$ of $\tau$, we can meaningfully talk about $\epsilon^\sigma_\tau|_{A_\bullet(\alpha)}$. Note that, under the above definition, we have the additivity relation
\begin{equation*}
\epsilon^\sigma_\alpha = \epsilon^\sigma_\tau|_{A_\bullet(\alpha)} + \epsilon^\tau_\alpha
\end{equation*} Indeed, at $\alpha$ this additivity relation follows by definition of $\epsilon^\sigma_\tau$. At a proper face $\beta$ of $\alpha$, the additivity relation follows by induction, using additivity in the lower-dimensional cases.

Moreover, whenever $\alpha$ is a proper face of $\tau$, the chain homotopy equation 
\begin{equation*}
\partial \epsilon^\sigma_\tau(\alpha) =\pi^\tau(\alpha)-\pi^\sigma(\alpha)-\epsilon^\sigma_\tau(\partial \alpha)
\end{equation*}
follows by induction, using additivity and the chain homotopy equations for $\epsilon^\sigma_\alpha$ and $\epsilon^\tau_\alpha$.

To complete the definition of $\epsilon^\sigma_\tau$, it remains to define $\epsilon^\sigma_\tau(\tau)$. If $\tau$ is a point, define $\epsilon^\sigma_\tau(\tau)$ to be the constant $1$-simplex at $\tau$. Now suppose $\tau$ is an $i$-face for $i>0$. We note that the chain $\pi^\tau(\tau)-\pi^\sigma(\tau)-\epsilon^\sigma_\tau(\partial \tau)$ is a cycle in $C_\bullet^\eta(|\sigma|)$, so by contractibility of $\sigma$ and step 2, we can let $\epsilon^\sigma_\tau(\tau)$ be a chain in $C_\bullet^\eta(|\sigma|)$ whose boundary is $\pi^\tau(\tau)-\pi^\sigma(\tau)-\epsilon^\sigma_\tau(\partial \tau)$. \qed
\\
\\Construction 3.2: we construct for each simplex $\sigma$ a chain homotopy $h^\sigma$ between $\pi^\sigma\colon A_\bullet(\sigma)\to C_\bullet(|\sigma|)$ and the inclusion, such that:

i) (Naturality on isomorphisms) Given an order-preserving linear isomorphism $f:\Delta^i\cong \sigma$, the homotopy $h^{\sigma}$ is given by conjugating $h^{\Delta^i}$ by the isomorphism $f$.

ii)For any proper face $\tau$ of $\sigma$, we have $h^\sigma|_{A_\bullet(\tau)}-h^\tau=\epsilon_\tau^\sigma$.

iii)For any face $\tau$ of $\sigma$ that lies in $C_\bullet^\eta(|\sigma|)$, we have $h^\sigma(\tau)\in C_\bullet^\eta(|\sigma|)$ as well.
\\
\\Proof: Condition i) serves as a definition in the case that $\sigma$ is not a standard simplex. So we may assume $\sigma=\Delta^k$.

If $\tau$ is a proper face of $\sigma$, we define $h^\sigma(\tau)=\epsilon^\sigma_\tau(\tau)+h^\tau(\tau)$. 

We have yet to define $h^\sigma(\sigma)$, but for any proper face $\tau$ of $\sigma$, we can meaningfully talk about $h^\sigma|_{A_\bullet(\tau)}$. We have
\begin{equation*}
h^\sigma|_{A_\bullet(\tau)}=\epsilon^\sigma_\tau+h^\tau
\end{equation*}
Indeed, at $\tau$ this follows by definition. At proper faces of $\tau$, this follows by using the additivity of the $\epsilon$ homotopies and using condition ii) for $h^\tau$ by induction.

If $\tau$ is a proper face of $\sigma$, the chain homotopy equation
\begin{equation*}
\partial h^\sigma(\tau)=\tau - \pi^\sigma(\tau) - h^\sigma(\partial \tau)
\end{equation*}
Follows by using condition ii) and induction.

We now check condition iii) for proper faces $\tau\subsetneq \sigma$. Indeed, if $\tau\in C_\bullet^\eta(|\sigma|)$, then by induction $h^\sigma(\tau)=\epsilon^\sigma_\tau+h^\tau(\tau)$ is a sum of chains in $C_\bullet^\eta(|\sigma|)$ which is therefore also in $C_\bullet^\eta(|\sigma|)$.

To complete the definition of $h^\sigma$, we define $h^\sigma(\sigma)$. If $\sigma=\Delta^0$, define $h^\sigma(\sigma)$ to be the constant $1$-simplex at $\sigma$. Now suppose $\sigma=\Delta^k$ for $k>0$. Note the chain $\sigma-\pi^\sigma(\sigma)-h^\sigma(\partial \sigma)$ is a cycle, so by contractibility of $\sigma$, we let $h^\sigma(\sigma)$ be a chain in $\sigma$ whose boundary is $\sigma-\pi^\sigma(\sigma)-h^\sigma(\partial \sigma)$. If $\sigma\in C_\bullet^\eta(|\sigma|)$, then we note that the chain $\sigma-\pi^\sigma(\sigma)-h^\sigma(\partial \sigma)$ in fact lies in $C_\bullet^\eta(|\sigma|)$, so by step 2, we can pick $h^\sigma(\sigma)$ to also lie in $C_\bullet^\eta(|\sigma|)$.\qed
\\
\\The above constructions only have naturality with respect to isomorphisms. The next construction will be natural with respect to all face maps between simplices.
\\
\\Construction 3.3: For each simplex $\Delta^k$ equippped with a nesting $\eta$, we construct a map of complexes $\rho_{\Delta^k,\eta}\colon A_\bullet(\Delta^k)\to C_\bullet^{\eta}(\Delta^k)$ and a chain homotopy $H_{\Delta^k,\eta}:A_n(\Delta^k)\to C_{n+1}(\Delta^k)$ between $i\circ \rho_{\Delta^k,\eta}$ and inclusion, such that:

i) (Naturality) For any face map $f\colon \Delta^i\hookrightarrow \Delta^k$, we have  $H_{\Delta^k,\eta}\circ f_* = f_*\circ H_{\Delta^i,f^*\eta}$.

ii) If $\sigma$ is a face of $\Delta^k$ that lies in $C_\bullet^\eta(\Delta^k)$, then $H_{\Delta^k,\eta}(\sigma)$ does as well.
\\
\\Proof: We will suppress the subscripts on $H$ and $\rho$. We first define $H$. Given a face $\sigma$ of $\Delta^k$, we define $H(\sigma)=h^\sigma(\sigma)\in C_\bullet(|\sigma|)\subset C_\bullet(\Delta^k)$ and define $\rho$ by $\rho(\sigma)=\sigma - H(\partial \sigma) - \partial H(\sigma)$. Naturality follows by naturality on isomorphisms in construction 3.2.

Note $\rho$ is indeed a map of complexes, since $\partial \rho(\sigma)=\partial \sigma - (\partial \sigma - \rho(\partial \sigma)) = \rho(\partial \sigma)$. The fact that $H$ is a chain homotopy between $i\circ \rho$ and inclusion is by definition.

It remains to check condition ii). Indeed, we can rewrite 
\begin{equation*}
\rho(\sigma) = \sigma - H(\partial \sigma) - (\sigma - \pi^{\sigma}(\sigma) - h^\sigma(\partial \sigma)) = \pi^\sigma(\sigma)+h^\sigma(\partial \sigma) - H(\partial \sigma)
\end{equation*}
But $\pi^\sigma(\sigma)\in C_\bullet^\eta(\Delta^k)$ and $h^\sigma(\partial \sigma)-H(\partial \sigma)\in C_\bullet^\eta(\Delta^k)$ since for each simplex $\tau$ in the boundary of $\sigma$, we have $h^\sigma(\tau)-H(\tau)=h^\sigma(\tau)-h^\tau(\tau)=\epsilon^\sigma_\tau(\tau)\in C_\bullet^\eta(\Delta^k)$. \qed
\\
\\Finally, we generalize to a topological space $X$ equipped with a nesting $\eta$. Let $\sigma\colon \Delta^k\to X$ be a simplex. Then we let $\rho_{\Delta^k}$ and $H_{\Delta^k}$ be given by construction 3.3, applied to $\Delta^k$ with the nesting $\sigma^*\eta$. We define $\rho(\sigma)=\sigma_*(\rho_{\Delta^k})$ and $H(\sigma)=\sigma_*(H_{\Delta^k}(\Delta^k))$. Since $\rho_{\Delta^k}$ has image in $C_\bullet^{\sigma^*\eta}(\Delta^k)$, we see that $\rho$ has image in $C_\bullet^\eta(\Delta^k)$. Similarly, whenever $\sigma\in C_\bullet^\eta(X)$, we have $H(\sigma)\in C_\bullet^\eta(X)$. Moreover, since each $\rho_{\Delta^k}$ is a map of complexes and each $H_{\Delta^k}$ is a chain homotopy between $i\circ \rho_{\Delta^k}$ and the inclusion, it follows by the naturality of $\rho_{(-)}$ and $H_{(-)}$ with respect to face maps that, similarly, $\rho\colon C_\bullet(X)\to C_\bullet^\eta(X)$ is a map of complexes and $H:C_n(X)\to C_{n+1}(X)$ is a chain homotopy between $i \circ \rho$ and the identity. On the other hand, the restriction of $H$ to $C_\bullet^\eta(X)$ gives a chain homotopy between $\rho\circ i$ and the identity on $C_\bullet^\eta(X)$, thus establishing that $i$ is indeed a chain homotopy equivalence.

\end{document}